\documentclass[12pt,a4paper,reqno]{amsart}
\usepackage{amssymb}
\usepackage{txfonts}
\usepackage{mathrsfs}
\usepackage{bbm}
\def\bc{\begin{center}}
\def\ec{\end{center}}
\def\be{\begin{equation}}
\def\ee{\end{equation}}

\def\N{\mathbb N}

\newtheorem{lem}{Lemma}[section]

\newtheorem{pro}[lem]{Proposition}
\newtheorem{thm}[lem]{Theorem}

\theoremstyle{remark}
\newtheorem{rem}{Remark}
\numberwithin{equation}{section}
\usepackage{colortbl}

\begin{document}
\title[Fast Khintchine spectrum ] {On the fast Khintchine spectrum in continued fractions}

\thanks{$^\dag$ Corresponding author.}

\author{Aihua Fan}
\address{Lamfa, Umr 7352 (Ex 6140), CNRS,
Universit\'e de Picardie Jules Verne, 33, Rue Saint Leu, 80039
Amiens Cedex 1, France} \email{ai-hua.fan@u-picardie.fr}

\author{Lingmin Liao}
\address{Lama, Umr 8050, CNRS, Universit\'e Paris-Est Cr\'eteil Val de Marne,
61, avenue du G\'en\'eral de Gaulle 94010 Cr\'eteil Cedex France
}\email{lingmin.liao@u-pec.fr}

\author{Baowei WANG$^\dag$}
\address{School of Mathematics and Statistics, Huazhong University of Science and Technology, 430074 Wuhan, China}
\email{bwei\_wang@yahoo.com.cn}

\author{Jun WU}
\address{School of Mathematics and Statistics, Huazhong University of Science and Technology, 430074 Wuhan, China}
\email{wujunyu@public.wh.hb.cn}
\begin{abstract}
For $x\in [0,1)$, let $x=[a_1(x), a_2(x),\cdots]$ be its continued
fraction expansion with partial quotients $\{a_n(x), n\ge 1\}$. Let $\psi : \mathbb{N} \rightarrow \mathbb{N}$ be a function with $\psi(n)/n\to \infty$ as $n\to
\infty$. In this note, the fast Khintchine spectrum, i.e., the Hausdorff dimension of the set
$$ E(\psi):=\Big\{x\in [0,1):
\lim_{n\to\infty}\frac{1}{\psi(n)}\sum_{j=1}^n\log a_j(x)=1\Big\}
$$ is completely determined without any
extra condition on $\psi$. 
\end{abstract}


 \maketitle

\addtocounter{section}{0}
\section{Introduction}
Continued fraction expansions are induced by the Gauss
transformation $T:[0,1)\to [0,1)$ given by  $$ T(0):=0, \
T(x)=\frac{1}{x} \ {\rm{(mod \ 1)}}, \ {\rm{for}}\ x\in (0,1).
$$ Let $a_1(x)= \lfloor x^{-1}\rfloor$ ($\lfloor \cdot \rfloor $ stands for the integral part)
  and $a_n(x)=a_1(T^{n-1}(x))$ for $n\ge 2$. Each irrational
number $x\in [0,1)$ admits a unique infinite continued fraction
expansion of the form
\begin{eqnarray}\label{ff1}
x=\frac{\displaystyle 1}{\displaystyle a_1(x)+ \frac{\displaystyle
1}{\displaystyle a_2(x)+\frac{\displaystyle 1}{\displaystyle
a_3(x)+\ddots}}}.
\end{eqnarray}
Sometimes,  (\ref{ff1}) is written as $x=[a_1,a_2,\cdots]$. The
integers $a_n$ are called the partial quotients of $x$.
 The $n$-th convergent $p_n(x)/q_n(x)$ of $x$
is given by $p_n(x)/q_n(x)=[a_1,\cdots, a_n]$.

The continued fraction is tightly connected with the classic
Diophantine approximation. For example, for any $v\ge 2$, the
well-known Jarn\'ik set
$$ \left\{x: \left|x-p/q\right|<q^{-v}, \ \ {\text{for infinitely many }}\ (p,q)\in
\mathbb{Z}^2\right\}
$$ is equal to $J_{v-2}$, where for any $\beta>0$, the set $J_{\beta}$ is defined by continued fractions as$$
J_{\beta}:=\left\{x: a_{n+1}(x)\ge q_n(x)^{\beta}, \ \ {\text{for
infinitely many }}\ n\in \mathbb{N}\right\}.
$$
%
%

\smallskip

The Gauss transformation is identified with an infinite symbolic
dynamical system if we consider the partial quotients as symbols.
The appearence of infinite symbols brings us new phenomena in
relative to the case of finite symbols. For example, consider the
set
$$ \left\{x\in [0,1): \mathbb{A}\Big\{\ {\frac{1}{n}}\sharp\{1\le
j\le n: a_j(x)=1\} \Big\}_{n\ge 1}= [0,1]\right\}
$$ where $\mathbb{A}(E)$ denotes the set of the accumulation points of a set $E$.
The Hausdorff dimension of this set is $1/2$ (see \cite{LMW}), while
in $b$-adic expansion a similar set is of Haudorff dimension 0 (see
\cite{Ol}). Another example is that the multifractal spectrum of the
level sets of the Khintchine constant
\begin{equation*} \Big\{x\in [0,1):\lim_{n\to\infty}\frac{1}{n}\sum_{j=1}^n\log
a_j(x)=\xi\Big\}
\end{equation*} is neither concave nor convex \cite{Fan}.
Because of the difference from the finite symbolic dynamical systems
and of the observed new phenomena, continued fractions attracted
much attention. One can find rich properties of the continued
fraction dynamical system in \cite{FLM,Fan, KJ, KS, KPW, MU, Ma, PW}
and related works therein.


\medskip

Let $\psi:\mathbb{N}\to \mathbb{N}$. Define $$
E(\psi)=\left\{x\in [0,1): \lim_{n\to\infty}\frac{\log
a_1(x)+\cdots+\log a_n(x)}{\psi(n)}=1\right\}.
$$ When $\psi(n)=\lambda n$ for some $\lambda>0$, the set $E(\psi)$ is a level set of the classic
Khintchine constant. Besides a detailed spectrum analysis of the classic Khintchine constant in \cite{Fan},
the authors also studied the fast Khintchine spectrum, i.e. the Hausdorff dimension of $E(\psi)$
when $\psi(n)/n\to \infty$ as $n\to\infty$. But the result for the latter case is incomplete.
Only under the conditions that $\lim_{n\to
\infty}\frac{\psi(n+1)}{\psi(n)}=b$ and
$\lim_{n\to\infty}(\psi(n)-\psi(n-1))=\infty$, the dimension of
$E(\psi)$ was given \cite{Fan}.
 In this note, we show that these extra conditions are unnecessary for determining
 the dimension of $E(\psi)$ in the case of fast Khintchine spectrum.

 Two functions $\psi$ and $\tilde{\psi}$ defined on
$\N$ are said to be {\em equivalent} if
$\frac{\psi(n)}{\tilde{\psi}(n)}\to 1$ as $n\to \infty$.

\begin{thm}\label{tE}  Let $\psi:\N\to \N$ with ${\psi(n)}/{n}\to \infty$ as $n\to \infty$.
If $\psi$ is equivalent to an increasing function, then $E(\psi)\ne
\emptyset$ and
$$ \dim_HE(\psi)=\frac{1}{1+b}, \ \ {\text{with}} \ b=\limsup_{n\to
\infty}\frac{\psi(n+1)}{\psi(n)}.
 $$ Otherwise,
 $E(\psi)=\emptyset$.
\end{thm}

\begin{rem}The method used in \cite{Fan} does not apply to general $\psi$.
This is explained in Section 3 below.\end{rem}

\begin{rem}The upper bound of $\dim_HE(\psi)$ is the difficult part of the proof of Theorem \ref{tE}. As a byproduct
of the proof, we get that for any $\beta>0$, the Hausdorff dimension
of the set
$$ J^*_{\beta}:=\Big\{x\in J_{\beta}: \lim_{n\to\infty}\frac{\log
q_n(x)}{n}=\infty\Big\}
$$ is $1/(2+\beta)$, i.e. one half of the dimension of the Jarn\'ik
set $J_{\beta}$. 
A detailed explanation is given at the end of this paper.
\end{rem}

\section{Preliminary}
This section is devoted to fixing some notation, recalling some
elementary properties enjoyed by continued fractions and citing some
technical lemmas in dimension estimation.

Throughout this paper, we use $\lfloor \cdot\rfloor$ to denote the
integral part of a real number, $|A|$ the diameter of a set
$A\subset \mathbb{R}$, $\mathcal{H}^s$ the $s$-dimensional Hausdorff
measure, and $\dim_{H}$ the Hausdorff dimension of a subset of $[0,1)$.

Recall that for any irrational number $x\in [0,1)$, $p_n(x)$ and
$q_n(x)$ are the numerator and denominator of  the $n$-th convergent
of $x$. 
It is known that $p_n=p_n(x)$ and $q_n=q_n(x)$ can be obtained recursively by the following
relations.
\begin{equation}\label{ff2.1}
\begin{split}
p_n=a_n(x) p_{n-1}+p_{n-2}, \ \ q_n=a_n(x)
q_{n-1}+q_{n-2}
\end{split}
\end{equation}with the conventions $p_0=q_{-1}=0$ and
$p_{-1}=q_0=1$. For each $n\ge 1$,
\begin{equation}\label{7}
 p_{n-1}q_n-p_nq_{n-1}=(-1)^n.
\end{equation}

 For any $n \geq 1$ and $(a_1,a_2,\cdots,a_n)\in
\mathbb{N}^n$, define $$ I_n(a_1, a_2, \cdots, a_n)=\big\{x\in
[0,1):\ a_1(x)=a_1, \cdots, a_n(x)=a_n\big\},
$$ which is the set of points beginning with $(a_1,\cdots, a_n)$ in their continued fraction expansions,
and is called a {\em cylinder} of order $n$.

Note that $p_n$ and $q_n$ are determined  by the first $n$ partial
quotients of $x$. So all points in $I_n(a_1,\cdots, a_n)$  determine
the same $p_n$ and $q_n$. Hence sometimes, we write
$p_n=p_n(a_1,\cdots, a_n)$ and $q_n=q_n(a_1,\cdots, a_n)$ to denote
$p_n(x)$ and $q_n(x)$ for $x\in I_n(a_1,\cdots, a_n)$.
\begin{pro}[\cite{Kh}]\label{p2.1}
For any $n\geq 1$ and $(a_1,\cdots, a_n)\in \mathbb{N}^n$, let $q_n$
be given recursively by (\ref{ff2.1}).
 The cylinder $I_n(a_1,\cdots,a_n)$ is an interval with the endpoints $p_n/q_n$ and $(p_n+p_{n-1})/(q_n+q_{n-1})$.
Then
\begin{equation}\label{5}
\frac{1}{2q_n^2}\le \Big|I_n(a_1,\cdots,
a_n)\Big|=\frac{1}{q_n(q_n+q_{n-1})}\le
\frac{1}{q_n^2}.\end{equation}
For each $n\ge 1$, $q_n(a_1, \cdots,
a_n)\ge 2^{(n-1)/2}$ and \begin{equation}\label{ff13}
\prod_{k=1}^na_k\le q_n(a_1,\cdots,a_n)\le   2^n \prod_{k=1}^na_k.
\end{equation}
\end{pro}

Now we  mention some known results concerning the dimension of sets
in continued fractions. Let $\{s_n\}_{n\geq 1}$ be a sequence of
integers and $\ell\geq 2$ be some fixed integer. Set$$
F(\{s_n\}_{n=1}^{\infty};\ell)=\big\{x\in [0,1): s_n\leq a_n(x)<\ell
s_n, \ {\rm{for \ all}}\ n\geq 1\big\}.
$$\begin{lem}[\cite{Fan}]\label{l3.1}
Under the assumption that $\frac{1}{n}\sum_{k=1}^n s_k\to \infty$ as
$n\to \infty$, one has\begin{eqnarray*}
\dim_HF(\{s_n\}_{n=1}^{\infty};\ell)=\liminf_{n\to \infty}\frac{\log
(s_1s_2\cdots s_n)}{2\log (s_1s_2\cdots s_n)+\log s_{n+1}}.
\end{eqnarray*}
\end{lem}

\begin{lem}[\cite{Fan}]\label{l2.7}
$$\dim_H\left\{x\in [0,1): \limsup_{n\to \infty}\frac{\log q_n(x)}{n}=\infty\right\}= \frac{1}{2}.$$
\end{lem}

\medskip

\section{Proof of Theorem \ref{tE}}

 Notice that $E(\psi)=E(\tilde{\psi})$ if $\psi$ and $\tilde{\psi}$ are equivalent.
We can assume that $\psi$ is increasing because of the following
simple lemma. \begin{lem}
 The set $E(\psi)\neq \emptyset$ if and only if $\psi$ is equivalent to an increasing function.
\end{lem}
\proof  If $E(\psi)$ is nonempty, take an $x_0\in E(\psi)$. Then put
$$\tilde{\psi}(n)=\Big \lfloor \log a_1(x_0)+\cdots+\log a_n(x_0)\Big \rfloor $$ for all
$n\geq 1$. Clearly $\tilde{\psi}$ is increasing. The functions
$\psi$ and $\tilde{\psi}$ are equivalent.

On the other hand, if $\psi$ is increasing, we have a point $x\in
E(\psi)$ such that for each $n\ge 1$ $$ a_n(x)=\lfloor
e^{\psi(n)-\psi(n-1)+1}\rfloor.
$$
\hfill $\Box$

Now we can proceed  the proof of Theorem  \ref{tE} with the
assumption that $\psi$ is increasing. 

$\bullet$ {\bf {\sc Lower bound.}} Apply Lemma \ref{l3.1} to
$s_n=\lfloor e^{\psi(n)-\psi(n-1)} \rfloor$ and $\ell=2$. Let
$$ F=\left\{x\in [0,1): \Big \lfloor e^{\psi(n)-\psi(n-1)}\Big \rfloor\leq
a_n(x)<2\Big \lfloor e^{\psi(n)-\psi(n-1)}\Big \rfloor, \ {\rm{for \
all}}\ n\geq 1\right\}
$$ which is subset of $E(\psi)$. We get immediately that $$
\dim_HE(\psi)\ge \frac{1}{1+b}.
$$

$\bullet$ {\bf {\sc Upper bound.}} This is the main part of the
proof.

Let us first recall the method used in \cite{Fan} under the extra
condition that $\lim_{n\to \infty}\frac{\psi(n+1)}{\psi(n)}=b\ge 1$.
Especially when $b>1$, we constructed a set containing $E(\psi)$ by
posing precise restrictions on each partial quotients, namely
\begin{equation}\label{12}\Big\{x\in [0,1): e^{L_n}\le a_{n}(x)\le
e^{M_n}, {\text{when}}\ n\gg 1\Big\},\end{equation} where (with a
small $\epsilon>0$)
$$ L_n=\frac{\psi(n)}{1+\epsilon}-\frac{\psi(n-1)}{1-\epsilon} \ \ \text{ and } \ \
M_n=\frac{\psi(n)}{1-\epsilon}-\frac{\psi(n-1)}{1+\epsilon}.
$$
By a standard covering argument, together with $\lim_{n\to
\infty}\frac{\psi(n+1)}{\psi(n)}=b$, we get the exact upper bound of
the dimension of $E(\psi)$. But as far as a general function $\psi$
is concerned,
 the above argument
fails. For example, take $$\psi(n)=(k+2)!, \ {\text{when}}\  k!\le
n<(k+1)!.$$ Then the set in (\ref{12})  reads as
\begin{eqnarray*}\left\{x\in [0,1): \left\{
                                  \begin{array}{ll}
                                    e^{c_1(k+2)!}\le a_n(x)\le e^{c_2(k+2)!}, & \hbox{when $n=k!$;} \\
                                    1\le a_n(x)\le e^{c_3 (k+2)!}, & \hbox{when $k!<n<(k+1)!$.}
                                  \end{array}
                                \right.
\right\}
\end{eqnarray*}
for suitably chosen constants $c_1,c_2,c_3$. According to Lemma
\ref{l3.1}, this set  has Hausdorff dimension $\ge 1/2$. However,
the dimension of $E(\psi)$ is equal to zero by Theorem \ref{tE}.

\medskip Now we are going to prove the upper bound of $\dim_HE(\psi)$ for a general function $\psi$.
Since $\psi$ is  increasing, we always have $b\ge 1$. We distinguish
two cases: $b=1$ and $b>1$.

 {\bf Case \ 1.} $b=1$. Lemma \ref{l2.7} serves for this case. According to the estimation (\ref{ff13}),
 since $\psi(n)/n\to\infty$ as $n\to \infty$, we have $$
 \lim_{n\to\infty}\frac{\log q_n(x)}{\psi(n)}=\lim_{n\to\infty}\frac{\log a_1(x)+\cdots+\log a_n(x)}{\psi(n)}.
 $$Thus  Lemma \ref{l2.7} gives us
$$\dim_HE(\psi)\leq \frac{1}{2}=\frac{1}{1+b}.$$

 {\bf Case \ 2}. $b>1$.
Fix an $\epsilon>0$. Choose a sequence of integers $\{n_k\}_{k=1}^{\infty}
\subset \mathbb{N}$ with $n_1$ large enough and for each $k\geq 1$
one has
\begin{equation}\label{ff3.3}\psi(n_k+1)\geq \psi(n_k)b(1-\epsilon), \ \ n_k\le \epsilon \psi(n_k).
\end{equation}
For each $N\ge 1$, let $$ E_N(\psi)=\Big\{x\in [0,1):
(1-\epsilon)<\frac{1}{\psi(n)}\sum_{j=1}^n\log a_j(x)<(1+\epsilon),
\ \forall \ n\ge N\Big\}.
$$ Then $$E(\psi)\subset \bigcup_{N\ge 1}E_N(\psi).$$
To estimate  the dimension of $E_N(\psi)$ for $N\ge 1$, we proceed
in three steps.
\medskip

 Step i. {\bf Find a cover of $E_N(\psi)$}.
For any $n\geq N$, set
\begin{equation}\label{ff3.4}
D_n(\epsilon)=\Big\{(a_1,\cdots,a_n)\in \mathbb{N}^n:
(1-\epsilon)<\frac{1}{\psi(n)}\sum_{j=1}^n\log
a_j<(1+\epsilon)\Big\}.
\end{equation} For every $(a_1,\cdots,a_n)\in D_n(\epsilon)$,  we define $$
{D}_{n+1}\Big(\epsilon; (a_1,\cdots,a_n)\Big)=\Big\{a_{n+1}\in
\mathbb{N}: (a_1,\cdots,a_n,a_{n+1})\in D_{n+1}(\epsilon)\Big\}.
$$

Clearly, by the definition of $D_n(\epsilon)$, we have
\begin{equation}\label{f2}
E_N(\psi)\subset \bigcap_{n=N}^{\infty}\ \mathfrak{D}_n(\epsilon), \
{\text{with}}\ \mathfrak{D}_n(\epsilon)=\bigcup_{(a_1,\cdots,
a_n)\in D_n(\epsilon)}I_n(a_1,\cdots,a_n).
\end{equation}
Now instead of considering the intersections in (\ref{f2}) from
$n=N$ until $n=\infty$, we only consider the intersection of two
consecutive terms. Namely, for any $n\ge N$,
\begin{align*}
E_N(\psi)&\subset \Big(\mathfrak{D}_n(\epsilon)\cap
\mathfrak{D}_{n+1}(\epsilon)\Big)=\bigcup_{(a_1,\cdots,a_n)\in
D_n(\epsilon)}J_n(a_1,\cdots,a_n),
\end{align*}
 where $$ J_n(a_1,\cdots,a_n)=\bigcup_{a_{n+1}\in
{D}_{n+1}(\epsilon; (a_1,\cdots,a_n))}I_{n+1}(a_1,\cdots,a_n,
a_{n+1}).
$$ 
Hence, for each $n\ge N$, we get a cover of $E_N(\psi)$:
\begin{equation}\label{3} \Big\{J_n(a_1,\cdots,a_n):
(a_1,\cdots,a_n)\in D_n(\epsilon)\Big\}.
\end{equation}
Thus the $s$-dimensional Hausdorff measure of $E_N(\psi)$ can be
 estimated as
\begin{eqnarray}\label{f1}
\mathcal{H}^s(E_N(\psi))\le \liminf_{n\to\infty}
\sum_{(a_1,\cdots,a_{n})\in D_{n}(\epsilon)}\Big|J_{n}(a_1,\cdots,a_{n})\Big|^s.
\end{eqnarray}

As we shall see, $J_n(a_1,\cdots,a_n)$ is a union of
cylinders of order $(n+1)$, say $I_{n+1}(a_1,\cdots, a_n, a_{n+1})$
with $a$ taking large values (Lemma \ref{l4}). Using this fact, the
length of $J_n(a_1,\cdots,a_n)$ will be well estimated.
\medskip

Step ii. {\bf Lengths of $J_n(a_1,\cdots,a_n)$}.\
 We begin with a fact on $D_{n+1}(\epsilon;a_1,\cdots,a_n)$.
\begin{lem} For each $(a_1,\cdots,a_n)\in D_n(\epsilon)$,
$${D}_{n+1}\Big(\epsilon; (a_1,\cdots,a_n)\Big)\ne \emptyset.$$  \end{lem}
\proof This follows from the following simple constructions.

(a) If $\sum_{j=1}^n \log a_j>(1-\epsilon)\psi(n+1)$, we choose $a_{n+1}=1$.

 (b) If $\sum_{j=1}^n \log a_j\le (1-\epsilon)\psi(n+1)$, we can choose $$a_{n+1}=\Big\lfloor{\frac{e^{\psi(n+1)}}{a_1 \cdots a_n}}\Big\rfloor.$$
\hfill $\Box$

Recall that the sequence of integers $\{n_k\}_{k\ge 1}$ is given in (\ref{ff3.3}).
\begin{lem}\label{l4}
For any $(a_1,\cdots, a_{n_k})\in D_{n_k}(\epsilon)$ and
$a_{n_k+1}\in D_{n_k+1}\big(\epsilon, (a_1,\cdots, a_{n_k})\big)$,
we have
\begin{equation}\label{ff3.6} \log a_{n_k+1}\geq
(1-\epsilon)\Big(\frac{b(1-\epsilon)^2}{1+\epsilon}-1\Big)\log
q_{n_k}=: \beta \log q_{n_k},
\end{equation}
\end{lem}
\proof
By the definitions of $D_n(\epsilon)$ and the first inequality in (\ref{ff3.3}),
 for any $(a_1,\cdots, a_{n_k})\in D_{n_k}(\epsilon)$
and $a_{n_k+1}\in D_{n_k+1}\big(\epsilon, (a_1,\cdots,
a_{n_k})\big)$, one has
\begin{align}\label{1}
\sum_{j=1}^{n_k+1}\log a_j\geq
\psi(n_k+1)(1-\epsilon)&\geq \psi(n_k)b(1-\epsilon)^2\nonumber\\
&\geq \frac{b(1-\epsilon)^2}{1+\epsilon}\sum_{j=1}^{n_k}\log a_j.
\end{align}

On the other hand, by
(\ref{ff13}) and the second inequality in (\ref{ff3.3}), we get
\begin{equation}\label{2} q_{n_k}(a_1,\cdots, a_{n_k})\le
2^{n_k}\prod_{j=1}^{n_k}a_j\le
\left(\prod_{j=1}^{n_k}a_j\right)^{\frac{1}{1-\epsilon}}.
\end{equation}
Combining (\ref{1}) and (\ref{2}), we obtain the desired result. \hfill $\Box$
\medskip

Now return back to the cover of $E_N(\psi)$ given in (\ref{3})
especially when  $n=n_k$. We estimate the length of
$J_{n_k}(a_1,\cdots,a_{n_k})$ for every $(a_1,\cdots,a_{n_k})\in
D_{n_k}(\epsilon)$. For $n=n_k$, by (\ref{ff3.6}) and Proposition
\ref{p2.1},  we have
\begin{align*}
\big|J_{n}(a_1,\cdots,a_{n})\big|&\leq \sum_{a:a\ge
q_{n}^{\beta}}\Big|\frac{a\cdot p_{n}+p_{n-1}}{a \cdot
q_{n}+q_{n-1}}- \frac{(a+1)p_{n}+p_{n-1}}{(a+1)q_{n}+q_{n-1}}\Big|.
\end{align*} By (\ref{7}), for all $a\in \N$,  the differences
appearing in the series have the same sign depending only the parity
of $n$. Thus the series is telescopic. Since
$\frac{(a+1)p_n+p_{n-1}}{(a+1)q_n+q_{n-1}}$ tends to $p_n/q_n$ as
$a\to \infty$, we get $$ \Big|J_n(a_1,\cdots,a_n)\Big|\le
\bigg|\frac{q_{n}^{\beta}p_{n}+p_{n-1}}{q_{n}^{\beta}q_{n}+q_{n-1}}-\frac{p_{n}}{q_{n}}\bigg|
=\frac{1}{(q_n^{\beta}q_n+q_{n-1})q_n}\le \frac{1}{q_{n}^{2+\beta}}.
$$Consider the liminf in (\ref{f1}) along the
subsequence $\{n_k\}_{k\ge 1}$, then we obtain
\begin{align}\label{6}
\mathcal{H}^s(E_N(\psi))\le \liminf_{k\to\infty}\sum_{(a_1,\cdots,a_{n_k})\in
D_{n_k}(\epsilon)}\left(\frac{1}{q_{n_k}}\right)^{s(2+\beta)}.
\end{align}

The last step is devoted to estimating the summation in ({\ref{6}})
under a suitable choice of $s$.

\medskip

Step iii.  {\bf Bernoulli measures}. A family of measures $\mu_{t}$
defined on cylinders is constructed firstly. For each $t>1$ and for
any $(a_1,\cdots,a_n)\in \mathbb{N}^n$, set
\begin{equation}\label{ff3.1} \mu_t(I_n(a_1,\cdots,a_n))=e^{-nP(t)-t\sum_{j=1}^n\log
a_j}, \end{equation}where $e^{P(t)}= \zeta(t)= \sum\limits_{k=
1}^{\infty} k^{-t}$. By Kolmogorov's consistency theorem, $\mu_t$
can be extended into a probability measure on $[0,1)$.

 Fix
$\epsilon>0$. By the assumption that $\lim_{n\to \infty}\psi(n)/n=
\infty$, one can choose some integer $N(\epsilon)\in \mathbb{N}$
such that for all $n\geq N(\epsilon)$,
\begin{equation}\label{ff3.2} nP\big(1+\frac{\epsilon}{2}\big)\leq \frac{\epsilon}{2}(1-\epsilon)\psi(n).
\end{equation}

 We claim that for each $n\ge N(\epsilon)$ and
$(a_1,\cdots,a_n)\in D_n$,
\begin{align}\label{4}{q_{n}}^{-(1+\epsilon)} \le \mu_{(1+\epsilon/2)}\big(I_{n}(a_1,\cdots,a_{n})\big).\end{align}
More precisely, for any $(a_1,\cdots,a_n)\in D_n$, by (\ref{ff3.4})
and (\ref{ff3.2}), we have
\begin{eqnarray}\label{f3}
\frac{\epsilon}{2}\sum_{j=1}^{n}\log a_j\ge  n
P(1+\frac{\epsilon}{2}).\end{eqnarray} Thus by (\ref{ff13}) and then
(\ref{f3}), we get
\begin{align*}{q_{n}}^{-(1+\epsilon)} \leq
e^{-(1+\epsilon)\sum_{j=1}^{n}\log a_j}&\leq
e^{-nP(1+\frac{\epsilon}{2})-(1+\frac{\epsilon}{2})\sum_{j=1}^{n}\log
a_j}.\end{align*}

 Choose $s=\frac{1+\epsilon}{2+\beta}$ in (\ref{6}).  By
(\ref{4}), we have
\begin{eqnarray*}
\mathcal{H}^{\frac{1+\epsilon}{2+\beta}}(E_N(\psi))\leq
\liminf_{k\to \infty}\sum_{(a_1,\cdots,a_{n_k})\in
D_{n_k}(\epsilon)}\mu_{(1+\epsilon/2)}\big(I_{n_k}(a_1,\cdots,a_{n_k})\big)\leq
1.
\end{eqnarray*}
Hence $$\dim_HE(\psi)\leq \sup_{N\ge 1}
\Big\{\dim_HE_N(\psi)\Big\}\leq \frac{1+\epsilon}{2+\beta}.$$
 Then the desired result
follows by letting $\epsilon\to 0$.
%
%
%
\hfill $\Box$

\smallskip
\noindent{\it Final remark}: Now we give a remark on the dimension
of $J^*_{\beta}$ and that of $J_{\beta}$. Recall that $J^*_{\beta}$
and $J_{\beta}$ are defined in Section 1. For any
$(a_1,\cdots,a_n)\in \N^n$, we define
$$
\tilde{J}_n(a_1,\cdots, a_n)=\bigcup_{a_{n+1}\ge
q_n^{\beta}}I_{n+1}(a_1,\cdots,a_n,a_{n+1}).
$$ Then it is clear that $$
J_{\beta}=\bigcap_{N=1}^{\infty}\ \bigcup_{n=N}^{\infty}\
\bigcup\tilde{J}_n(a_2,\cdots,a_n),
$$ where the last union is taken over all
$(a_1,\cdots,a_n)\in \N^n.$ While $$
J^*_{\beta}\subset\bigcap_{N=1}^{\infty}\ \bigcup_{n=N}^{\infty}\
\bigcup\tilde{J}_n(a_2,\cdots,a_n),
$$ where the last union is taken over all
$(a_1,\cdots,a_n)\in \N^n$ with  $\frac{\log q_n}{n}$ being
sufficiently large.
 As a result,
\begin{align*}
    \mathcal{H}^s(J_{\beta})\le\liminf_{N\to\infty}\sum_{n=N}^{\infty}\ \ \sum_{(a_1,\cdots,a_n)\in
    \N^n}\left(\frac{1}{q_n}\right)^{s(2+\beta)},
\end{align*}while $$   \mathcal{H}^s(J^*_{\beta})\le\liminf_{N\to\infty}\ \sum_{n=N}^{\infty}\ \ \sum_{(a_1,\cdots,a_n)\in
    \N^n, (\log q_n)/n \
    {\text{large}}}\left(\frac{1}{q_n}\right)^{s(2+\beta)}.
$$
By (\ref{5}), we know that
\begin{equation}\label{10}
1\le \sum_{a_1, \cdots, a_n\in \N^n} {q_n}^{-2}\le 2.
\end{equation}While, by (\ref{4}), we get \begin{equation}\label{11}
\sum_{a_1,\cdots,a_n: \log q_n/n \ {\text{large}} }
{q_n}^{-(1+\epsilon)}\le 1.
\end{equation}
 Comparing of (\ref{10}) and (\ref{11}) reveals
that $$ \dim_HJ_{\beta}\le \frac{2}{2+\beta}, \ \ \
\dim_HJ^*_{\beta}\le \frac{1}{2+\beta}.$$

Actually we have proven that $\dim_HJ^*_{\beta}=\frac{1}{2+\beta}$
since $E(\psi)$ can serve as a subset of $J^*_{\beta}$.

 \vskip 10pt \noindent {\sc Acknowledgement:} This work was partially supported by PICS program No. 5727, RFDP20090141120007, NSFC 10901066 and NSFC 11171124. The authors thank the Morningside Center of Mathematics, Beijing for its hospitality.

 {\small

\end{document}